\documentclass[aip,reprint]{revtex4-1}
\usepackage{graphicx}
\usepackage{amsmath}
\usepackage{amsfonts}
\usepackage{amssymb}
\usepackage{hyperref}
\usepackage{algorithm}
\usepackage{algpseudocode}
\usepackage{mathtools}
\usepackage{tikz}      
\usepackage{float}

\newcommand\Algphase[1]{%
\vspace*{-.3\baselineskip}\Statex\hspace*{\dimexpr-\algorithmicindent-2pt\relax}
\Statex\hspace*{-\algorithmicindent}\textbf{#1}%
\vspace*{-.7\baselineskip}\Statex\hspace*{\dimexpr-\algorithmicindent-2pt\relax}\rule{0.48\textwidth}{0.4pt}
}

\begin{document}
\title{Parallel replica dynamics method for bistable stochastic reaction networks: simulation and sensitivity analysis}
\author{Ting Wang}
\email{tingw@udel.edu}
\author{Petr Plech\'{a}\v{c}}
\email{plechac@udel.edu}
\affiliation{Department of Mathematical Sciences, University of Delaware, Delaware 19716 USA}

\begin{abstract}
Stochastic reaction networks that exhibit bi-stable behavior are common in many fields such as systems biology and materials science.
Sampling of the stationary distribution is crucial for understanding and characterizing the long term dynamics of bistable stochastic dynamical systems.
However, this is normally hindered by the insufficient sampling of the rare transitions between the two metastable regions.
In this paper, we apply the parallel replica (ParRep) method for continuous time Markov chain \cite{ParRep-CTMC} to accelerate the stationary distribution sampling of bistable stochastic reaction networks. 
The proposed method uses parallel computing to accelerate the sampling of rare transitions and
it is very easy to implement.  
We combine ParRep with the path space information bounds \cite{dupuis2016path} for parametric sensitivity analysis. 
We demonstrate the efficiency and accuracy of the method by studying the Schl\"{o}gl model and the genetic switch network.  
\end{abstract}

\maketitle



\section{Introduction}
Stochastic reaction networks have become increasingly important as a tool for modeling complex biological and chemical systems with random noises \cite{mcadams1999sa}.
Simulation of real-world reaction networks using the stochastic simulation algorithm (SSA) \cite{gillespie-SSA} can be computationally intractable due to the multiscale feature of the systems.
For instance, reaction networks in biological cells often involve vastly different numbers of molecules for different species and rate constants for different reaction channels\cite{kang2013separation}.
Therefore, the system is metastable in the sense that the SSA rarely samples the reactions 
involving small rate constants or low population species.  
Our paper addresses with another type of metastable issue associated with reaction networks. 
We consider metastablity that is caused by extremely 
rare transitions between two separate regions of the state space, i.e., 
the bistable reaction networks. 
It has been discovered recently that many biological and physical systems exhibit bistability and hence it is of great interests to understand the bistable phenomenon. \cite{mehta2008exponential, gardner2000construction, umulis2006robust, angeli2004detection}

We study two important aspects regarding bistable reaction networks: accelerated stationary distribution sampling and parametric uncertainty quantification, \cite{dupuis2016path, pantazis2013relative, arampatzis2015accelerated} using the parallel replica dynamics (ParRep) method. \cite{voter1998parallel, ParRep-SDE, ParRep-Chain, ParRep-CTMC}
We know that SSA based sampling for bistable reaction networks can be extremely expensive  
because of the rare sampling of transitions between two metastable regions.
As a remedy for this issue, the ParRep uses multiple parallel replicas to explore the transition path between the two metastable regions with a controllable error.  
The method was originally designed for sampling rare events in molecular dynamics simulation  such as Langevin dynamics. \cite{voter1998parallel}
The mathematical framework of ParRep was recently developed for the discrete time Markov chains (DTMC).\cite{ParRep-SDE, binder2015generalized, ParRep-Chain, aristoff2015parallel}
In this paper we apply the version of ParRep algorithm that we developed for continuous time Markov chains (CTMC) \cite{ParRep-CTMC} to accelerate the simulation of bistable reaction networks. Furthermore, the algorithm allows us to efficiently sample the stationary distribution starting from the transient regime. 
We also investigate the parametric sensitivity problem of bistable reaction networks. 
Basically, we study the change of bistable system outputs to perturbations in system parameters. 
This enables us to quantify the parametric uncertainty and system robustness.  
We point out that the proposed version of ParRep can be easily combined with the path space information bounds \cite{dupuis2016path} to provide useful information and reductions for parametric sensitivity analysis in high dimensions.

\subsection{Stochastic reaction network model}\label{sec:SRN}
We consider a well-mixed chemical system with $n$ species interacting through $m$ reaction channels with system size $V$.
Under the well-mixed assumption, the molecular population is modeled as 
an $n$ dimensional CTMC $X^V(t)$.
The numbers of molecule of the $i$th species consumed and produced in the $j$th reaction are denoted by $\eta_{ij}^-$ and $\eta_{ij}^+$, respectively.
We call the net change $\eta_j = \eta_j^- - \eta_j^+$ caused by the $j$th reaction the stoichiometric vector, which is independent of the system size $V$.
Each reaction channel is associated with a propensity function $\lambda_j^V(x, c), j = 1, \cdots, m$ such that given $X^V(t) = x$, the probability of the $j$th reaction occurs at the infinitesimal time interval $[t, t + \delta t)$ is $\lambda_j^V(x, c) \delta t$, where $c$ is the vector of rate constants in $\mathbb{R}^l$.
{\it In this paper, we will suppress $c$ when we write the propensity functions unless
we study the parametric sensitivity with respect to $c$.}
From the propensity functions, we can construct the transition rate matrix (or the infinitesimal generator) $Q^V$ of the Markov chain $X^V$ such that 
\begin{equation}\label{eqn:Q}
Q_{x, y}^V = \begin{cases}
\lambda_j^V(x,c) & y = x + \eta_j~\text{for some}~j = 1, \ldots, m;\\
0                        & \text{Otherwise.}
\end{cases}
\end{equation}
Moreover, it is well known that the time evolution of $X^V$ is characterized by the random time change representation\cite{ethier2009markov}
\begin{equation}
X^V(t) = X^V(0) + \sum_{j = 1}^m \mathcal{P}_j\left(\int_0^t \lambda_j^V(X^V(s)) \,ds\right) \eta_j,
\end{equation}
where $\mathcal{P}_j$ are independent unit rate Poisson processes.

For a fixed system of the size $V$, the probability distribution of the population process $X^V$ is completely governed by the chemical master equation (CME)
\begin{equation}\label{eqn:CME}
\frac{dp^V(x, t)}{dt} = \sum_{j=1}^m \lambda_j^V(x - \eta_j) p^V(x - \eta_j, t) 
- \lambda_j^V(x) p^V(x, t),
\end{equation}
where $p^V(x, t) = \mathbb{P}(X^V(t) =x)$.
In principle, the CME enables the computation of the distribution of $X^V(t)$ for any $V$. However, the CME is normally an infinite dimensional system which cannot be solved explicitly in general.
Therefore, Monte Carlo methods such as Gillespie's SSA are commonly used to obtain the numerical solution to the CME. 

We denote by $X_V(t) = V^{-1}X^V(t)$ the corresponding concentration process for a system with size $V$.  
When $V$ is large, the randomness of the reaction network can be neglected and $X_V(t)$ can be approximated by the solution of the reaction rate equation (RRE) \cite{kurtz1970solutions}
\begin{equation}\label{eqn:RRE}
\frac{d\bar{x}}{dt} = \sum_{j=1}^m \eta_j \lambda_j(\bar{x})
\end{equation}
in the sense that
\begin{equation}\label{eqn:fluid-limit}
\lim_{V\to \infty}\sup_{0 \leq s \leq t} |X_V(s) - \bar{x}(s)| = 0
\end{equation}
almost surely for any $t > 0$, 
where we assume $\lambda_j(x) = V^{-1} \lambda_j^V(V x)$ for all $j = 1, \ldots, m$. 
Throughout this paper, we assume such form for propensity functions and hence the random time change representation for the concentration process is
\begin{equation*}
X_V(t) = X_V(0) + \frac{1}{V}\sum_{j = 1}^m \mathcal{P}_j\left(\int_0^t V\lambda_j(X_V(s)) \,ds\right) \eta_j.
\end{equation*}

We focus on reaction networks that are modeled by an ergodic CTMC $X^V(t)$ such that the stationary distribution $\pi^c$ exists and the ergodic limit
\[
\lim_{t\to\infty}\frac{1}{t}\int_0^t f(X^V(s)) \,ds = \pi^c(f)
\]
holds for suitable observables $f$.
Here the stationary distribution $\pi_V^c$ depends on $c$ since the process $X^V(t)$ depends on $c$.
The gradient of $\pi_V^c(f)$ with respect to the parameter $c$, i.e., $\nabla \pi_V^c(f) $, serves as an indicator for the system's parametric uncertainty or robustness. 
We call the estimation of 
\[\nabla \pi_V^c(f) = \left(\frac{\partial \pi_V^c(f)}{\partial c_1}, \ldots, \frac{\partial \pi_V^c(f)}{\partial c_l}\right)^{\text{tr}}\] 
the stationary sensitivity analysis problem, where the superscript ``tr" means transpose.

\subsection{Reaction networks with bistability}
In this paper, we are mainly interested in accelerating simulation and sensitivity analysis for bistable reaction networks, i.e., 
reaction networks whose RRE has a pair of asymptotically stable fixed points $\bar{x}_+$ and $\bar{x}_-$ separated by a saddle point $\bar{x}_0$. 
We denote the neighborhood of $\bar{x}_+$ by $W_+$ and the neighborhood of $\bar{x}_-$ by $W_-$. 
If we neglect the randomness of the network, any initial point that is placed in $W_+$ (resp. $W_-$) approaches to $\bar{x}_+$ (resp. $\bar{x}_-$) eventually. 
However, due to the random noise (since $V$ is finite), the system is subject to rare, large fluctuations which make the concentration process $X_V(t)$ to be far away from one stable fixed point and enter into the neighborhood of the other stable point.
The theoretical tool to study this type of large fluctuations is the large deviation principle (LDP) \cite{LDP-Dembo-Zeitouni, LDP-Shwartz-Weiss, LDP-Dupuis-Ellis}.
The key ingredient in the LDP of $X_V$ is the rate function (or action) $I_0^T(x)$ which characterizes the exponentially small probability for $X_V$ remaining in a small neighborhood of a path $x$, i.e.,
\[\mathbb{P}\left(\sup_{0\leq t \leq T} |X_V(t) - x(t)| \geq \delta\right) \approx e^{-VI_0^T(x)}\]
for all small $\delta$ when $V$ is large. 
By minimizing the rate function over the path space one can find the so called most probable  path \cite{dykman1994large}.
In a bistable reaction network, $X_V$ sojourns in $W_+$ (resp. $W_-$) for long time until there is an exponentially small probability for it to leave $W_+$ (resp. $W_-$) along the most probable path.
In this sense, we call $W_+$ and $W_-$ metastable sets for $X_V$ since the sojourn times in both sets are exponentially long.
The metastability issue normally leads to insufficient sampling of transition events between $W_+$ and $W_-$ and consequently makes it computationally prohibitive to sample the stationary distribution $\pi_V^c$ for $X^V$.    
In this work we aim to speed up the sampling of $\pi_V^c$ by accelerating the exit from metastable sets using parallel computing.

\section{methodology}\label{sec:ParRep}
\subsection{Parallel replica dynamics}
The idea of ParRep was first introduced for simulating rare events \cite{voter1998parallel} and was recently formalized in several papers\cite{ParRep-SDE, ParRep-Chain, ParRep-CTMC}. 
Our goal in this section is to introduce the ParRep method \cite{ParRep-CTMC} to
accelerate the simulation of bistable stochastic reaction networks and estimate the stationary distribution. 
{\it Since we are considering fixed volume $V$ in this section, we will suppress the superscript $V$ to simplify the notations}. 

The theoretical justification for ParRep relies on the notion of the quasi-stationary distribution (QSD). 
Given a set $W$ and a DTMC $X_n$, a distribution $\nu$ is
called the quasi-stationary distribution of $X_n$ in $W$ if 
\begin{equation}
\nu(A) = \mathbb{P}^{\nu}(X_n \in A | N > n)
\end{equation} 
for all $n = 1,2,\ldots$ and any measurable set $A \subset W$, where $N$ is the first exit time of $X_n$ from $W$. 
The definition roughly says that the QSD is a distribution supported on $W$ such that if the initial distribution is $\nu$, then the DTMC $X_n$ remains distributed with $\nu$ before it exits $W$.
The existence and uniqueness of the QSD in this setting can be shown rigorously \cite{ParRep-CTMC}.
The consequence of assuming that $X_n$ starts at the QSD $\nu$ in $W$ is that 
the first exit time $N$ follows a geometric distribution with some parameter $1-\lambda$, i.e., 
$\mathbb{P}^{\nu}(N > n) = \lambda^n$ for all $n = 1, 2, \ldots$.
Moreover,  the first exit time $N$ and the exit state $X_N$ are independent\cite{QSD}. 

Now suppose we have $R$ independent and identically distributed replicas $(X_n^1, \ldots, X_n^R)$ of $X_n$, each with initial distribution QSD $\nu$.
Denote the first exit time of the $r$th replica by $N^r$ and define the smallest first exit time among all $R$ replicas by 
\[N^* = \min_{r} N^r.\]
Note that there could be more than one replicas which realize the $N^*$ (exit after the same number of steps), we denote by $K$ the smallest index among the exited replicas, i.e.,
\[K = \min\{r = 1,\ldots, R; X_{N^*}^r \notin W\}.\] 
Assuming each of the replicas of $X_n$ is initially distributed with the QSD $\nu$,
the following two results are crucial for the design of the ParRep algorithm\cite{ParRep-Chain}.
\begin{enumerate}
\item $X_{N^*}^K$ is independent of $R(N^* - 1) + K$;
\item $(X_{N^*}^K, R(N^* - 1) + K)$ has the same distribution as $(X_{N^1}^1, N^1)$.
\end{enumerate} 
The first result states that the first exit state from $W$ over $R$ replicas is independent with the total sojourn time over $R$ replicas. 
Furthermore, the second result guarantees that joint distribution of the first exit time and the first exit state is independent of the number of replicas.  
These facts suggest that we can use multiple replicas to explore a metastable region in order to accelerate the sampling of exit events but without changing the exit distribution.
That is, we can achieve acceleration by using parallel computing.
However, the gain of efficiency in this procedure is under the assumption 
that all replicas start with the QSD of $W$, which is not the case in general.
In order to sample the QSD for launching the parallel step, some preparation steps are needed to
make the process to be well into the quasi-stationary state.
Therefore, a complete cycle of ParRep can be roughly divided into three steps, 
\begin{itemize}
\item[S1] Decorrelation: simulate $X_n$ until the QSD $\nu$ of the current metastable set $W$ is sampled. Proceed to the dephasing step;
\item[S2] Dephasing: prepare a sequence of iid initial state $(x_1, \ldots, x_R)$ from $\nu$. 
Proceed to the parallel step;
\item[S3] Parallel: launch $R$ replicas of $X_n$ at $(x_1, \ldots, x_R)$ to explore the exit path from $W$. Return to the decorrelation step.
\end{itemize} 

We can adapt the above ParRep procedure for DTMC to the simulation of CTMC through simulating its embedded chain. 
More significantly, the algorithm can be modified to effectively sample the stationary distribution of a CTMC without the detailed balance assumption.   
We present the ParRep algorithm for CTMC in Algorithm \ref{alg:ParRep}.
The setup of notations in the ParRep algorithm is as follows. 
\begin{itemize}
\item[] $\tilde{X}(t)$: ParRep process we simulate throughout the ParRep algorithm;
\item[] $T_s$: time clock throughout the ParRep algorithm;
\item[] $I_s$: accumulated contribution to the time integral $\int_0^{T_{s}} f(\tilde{X}(s))\,ds$
throughout the ParRep algorithm;
\item[] $N_c$: count of transitions in each decorrelation step;
\item[] $n_c$: decorrelation threshold;
\item[] $N_p$: count of transitions in each dephasing step;
\item[] $n_p$: dephasing threshold;
\item[] $\tau$: holding time for the next reaction;
\item[] $J$: index of the next reaction;
\end{itemize}
Before we start the ParRep algorithm, we
choose fixed decorrelation threshold $n_c$ and dephasing 
threshold $n_p$ and initialize $T_s = 0, I_s = 0$ and $\tilde{X} = x_0$.

The procedure of the decorrelation step can be summarized as follows.
If $W$ is not a metastable set, then the process would 
leave $W$ rapidly and hence there is no need to launch the following dephasing and parallel steps. 
However, if $W$ is metastable then the process would remain in $W$ for at least $n_c$ transitions and the algorithm proceed to the dephasing step.
Since we assume $n_c$ is large enough for the process to reach the QSD of $W$, 
the state we obtain after $n_c$ transitions is asymptotically distributed according to the QSD.  
The dynamics in the decorrelation step is exact and hence there is no loss of accuracy and no acceleration either during this step. 
In the dephasing step, we apply the Fleming-Viot particle technique \cite{binder2015generalized} to sample 
a sequence of iid initial states that can be used in the subsequent parallel step.
Similar to the decorrelation step, we specify the dephasing threshold $n_p$ and let all $R$ replicas to evolve for $n_p$ transitions (jumps). 
During this procedure, if a replica leaves $W$ then we force it to restart from
the current state of another replica (chosen uniformly).
Similar to $n_c$, $n_p$ is large enough so that we sample a sequence of QSD distributed states
$(x_1,\ldots, x_R)$. 
Note that the dephasing step does not contribute anything to the $T_s$, $I_s$ and $\tilde{X}$,
its only purpose is to prepare the initial states $(x_1,\ldots, x_R)$ for the subsequent parallel step.
\begin{algorithm}[H]
  \caption{Parallel replica algorithm}
  \label{alg:ParRep}
   \begin{algorithmic}[1]
   \Algphase{Decorrelation Step:}
   \State initialize $N_c = 0$
   \While{$N_c < n_c$}
       \State generate $\tau$ and $J$
       \State update $T_s \gets T_s + \tau$, $I_s \gets I_s + f(\tilde{X}) \tau$, $\tilde{X} \gets \tilde{X} + \eta_J$
       \If{$\tilde{X}$ is still in $W$}
           \State $N_c \gets N_c + 1$
       \Else
           \State $N_c \gets 0$
       \EndIf
   \EndWhile
   \State proceed to dephasing step
   \Algphase{Dephasing Step:}
   \State initialize $N_p = 0$
   \State launch $R$ replicas of DTMC $X_n^1,\ldots, X_n^R$ from any initial distribution
   \While{$N_p < n_p$}
       \State generate $J^r$ for $r=1,\ldots, R$
       \State $X^r \gets X^r + \eta_{J^r}$, $N_p \gets N_p + 1$
       \If{any replica leaves $W$}
           \State randomly choose the state of a replica that remains
           \State  in $W$ and restart the exited replica from that state
       \EndIf
   \EndWhile
   \State set $(x_1, \ldots, x_R) = (X^1, \ldots, X^R)$
   \State proceed to parallel step with initial states $(x_1, \ldots, x_R)$    
   \Algphase{Parallel Step:}
   \State initialize $N^* = 0$, $K = R$
   \State launch $R$ replicas of CTMC $(X^1(t), \ldots, X^R(t))$ with initial states $(x_1, \ldots, x_R)$;
   \While{all $X^r$ in $W$}
       \State generate $\tau^r$ and $J^r$ for $r=1,\ldots, R$
       \State $N^* \gets N^* + 1$ 
       \If{any replicas leave $W$}
           \State $K \gets \min\{r = 1,\ldots, R; X_{N^*}^r \notin W\} $
       \EndIf
       \State $T_s \gets T_s + \sum_{r=1}^K \tau^r$
       \State $I_s \gets I_s + \sum_{r=1}^K f(X^r) \tau^r$
       \State $X^r \gets X^r + \eta_{J^r}$ for $r = 1, \dots, K$
   \EndWhile
   \State $\tilde{X} \gets X_{N^*}^K$
   \State return to the decorrelation step
    \end{algorithmic}
\end{algorithm}
The acceleration of ParRep comes from the parallel step. 
We launch $R$ parallel replicas from $(x_1, \ldots, x_R)$ to explore the exit event from $W$, that is, sample $N^*$, $K$ and the first exit state $X_{N^*}^K$.
Since $(X_{N^*}^K, R(N^* - 1) + K)$ has the same distribution as $(X_{N^1}^1, N^1)$, 
sampling of exit events with $R$ replicas (i.e., sample $N^*$ and $K$) in the parallel step is approximately $R$ times faster than that with serial simulation (i.e., sample $N^1$). 
Moreover, all the generated data from each replica in the parallel step are collected in order to 
sample the stationary distribution $\pi^c$.
This is through the update of the clock time $T_s$ and the time integral $I_s$.
Note that sampling $\pi^c$ by reusing these generated data from ParRep is statistically correct (asymptotically when $n_c$ and $n_p$ are large) comparing to that from the serial simulation.  
In fact, we have shown that the averaged contribution to $T_s$ or $I_s$ 
over each ParRep cycle is independent with the replica number $R$, provided that $x_1, \ldots, x_R$ are independent and distributed according to the QSD \cite{ParRep-CTMC, aristoff2015parallel}.

\subsection{Accuracy and efficiency}\label{sec:ParRep-accuracy-efficiency}
The accuracy of ParRep method relies on the choice of the decorrelation step $n_c$ and the dephasing step $n_p$ since these parameters determine how ``good" we sample the QSD before the parallel step.
In practice, we would never have exact sampling of the QSD at each ParRep cycle and hence there is an error associated with the inexact sampling of QSD.
However, for large $n_c$ and $n_p$ we can expect that the error is sufficiently small.
In fact, this can be justified by the following result\cite{ParRep-CTMC}.  
For fixed $n_c$, we define the distribution $\nu_{n_c}$ as 
\[\nu_{n_c}(A) = \mathbb{P}(X_{n_c} \in A|N>n_c)\]
for any measurable set $A \subset W$,
i.e., the distribution of $X_{n_c}$ conditioned on no exit event occurred after $n_c$ transition steps.
If we assume that the dephasing step is exact (i.e., $(x_1, \ldots, x_R)$ are independent and distributed as the QSD), then the averaged error for sampling $I_s$ over each ParRep cycle can be bounded by a constant times the total variation $\|\mu_{n_c} - \mu\|_{\text{TV}}$. 
Furthermore, the total variation converges geometrically fast in terms of $n_c$. 
This justifies that the dynamics of transition from one metastable set to another metastable set (i.e., one ParRep cycle) is asymptotically correct.
The analysis of the global error from all ParRep cycles is hard to analyze. 
However,  our numerical experiments in Sec. \ref{sec:numerical_eg} suggest that ParRep is a 
rather accurate algorithm for long time simulation. 

We briefly discuss the efficiency of ParRep for CTMC. 
In this paper, we define the speedup as the ratio between the total computational time of serial simulation and that of ParRep simulation. 
In the idealized scenario, the speedup factor of ParRep could be up to the number of replica used in the simulation as suggested by the properties. 
However, in practice the preparation for a sequence of QSD initial states offsets this linear acceleration.
Heuristically, the efficiency of ParRep relies on the metastability of the set.
If the set $W$ is strongly metastable, then the time spent in the decorrelation and dephasing steps is negligible comparing to the acceleration achieved in the parallel step. 
However, if the set $W$ is not truly metastable then the parallel step would not be activated and hence the ParRep is equivalent to SSA. 
In fact, this argument can be formalized and it turns out that the efficiency of ParRep 
is determined by the ratio $\lambda_1 / (\text{Re}(\lambda_2) - \lambda_1)$, where 
$\lambda_1$ and $\lambda_2$ (with $0>\lambda_1>\text{Re}(\lambda_2)$ ) are the two largest 
eigenvalues of the transition rate matrix $Q$ (see \eqref{eqn:Q} for definition) restricted to the metastable set $W$.
We do not pursue this aspect rigorously in this paper.
Interested reader could refer to the related literature \cite{binder2015generalized}.

\section{Path space information bound}
In this section, we combine the ParRep method with the path-space information bounds \cite{dupuis2016path} to accelerate the parametric sensitivity analysis of stochastic reaction networks. The bounds are derived using several concepts in information theory. 
For the readability of the paper, we briefly review these concepts and their connections in Appendix \ref{app:information-theory}.

Recall that we define the sensitivity analysis problem at the end of Section \ref{sec:SRN}.
There exist several types of sensitivity analysis methods such as the finite difference \cite{rathinam2010efficient, anderson2012efficient}, likelihood ration \cite{plyasunov2007efficient} and infinitesimal perturbation analysis or pathwise derivative method \cite{sheppard2012pathwise}.
We refer them as the direct methods since they aim to estimate the sensitivity directly.
However, direct estimation of the sensitivity can be extremely expensive due to their large variances \cite{wang2016efficiency} and complexity when applied to large reactions networks. 
Alternatively, we aim to compute a gradient-free upper bound of the sensitivity.
The computed sensitivity bounds can be used for screening out those insensitive parameters (with small bounds) and then direct methods can be applied for the remaining  of parameters \cite{arampatzis2015accelerated}.  

In general, given a probability distribution $P^c$ which depends on a vector of parameters $c$,   
we define the sensitivity index of an observable $f$ (along the direction $v$) as
\begin{equation}\label{eqn:SI}
S_{f, v}(P^{c}) = \left|\lim_{\epsilon \to 0}\frac{1}{\epsilon}\left\{\mathbb{E}_{P^{c+\epsilon v}}[f] - \mathbb{E}_{P^{c}}[f]\right\}\right|.
\end{equation}
Note that in the case that $P^c = \pi^c$ and $v = e_k$ (the $k$-th basis vector), the sensitivity 
index is simply the $k$-th component of the gradient $\nabla \pi^c(f)$.   
When we are interested in the sensitivity analysis of the stochastic process $X(t, c)$ with stationary distribution $\pi^c$, it is often convenient to interpret the distribution $P^c$ as the path space distribution $P_{[0,T]}^c$, i.e., the probability distribution of paths of $X(t, c)$ on the time interval $[0, T]$.
It can be shown that in the transient regime (i.e., the initial distribution of $X(t, c)$ is not $\pi^c$) the sensitivity index can be bounded by
\begin{equation}\label{eqn:transient-bounds}
S_{f, v}(P_{[0,T]}^{c}) \leq  \sqrt{\text{Var}_{P_{[0,T]}^{c}}(f)} \sqrt{v^{\text{tr}} \mathcal{I}(P_{[0, T]}^c) v},
\end{equation}
where $\mathcal{I}(P_{[0, T]}^c)$ is the path space Fisher information matrix (FIM) of the relative entropy $\mathcal{R}(P_{[0,T]}^c | P_{[0,T]}^{c+\epsilon v})$ (see Appendix \ref{app:UQ-bounds} for a formal derivation). 
In the stationary regime, a similar sensitivity bound can be derived.
That is, the stationary sensitivity index can be bounded by 
\begin{equation}\label{eqn:stationary-bounds}
S_{f,v}(\pi^{c}) \leq \sqrt{\tau_{\pi^c}(f)} \sqrt{v^{\text{tr}} \mathcal{I}_\mathcal{H}(P^{c}) v},
\end{equation} 
where $\tau_{\pi^c}(f)$ is the integrated auto-correlation function (IAF)   
and $\mathcal{I}_\mathcal{H}(P^{c})$ is the path space FIM of the relative entropy rate 
$\mathcal{H}(P^c | P^{c+\epsilon v})$. 
In fact, $\mathcal{I}_\mathcal{H}(P^{c})$ can be roughly interpreted as $\lim_{T\to\infty} T^{-1} \mathcal{I}(P_{[0, T]}^c)$.
See Appendix \ref{app:UQ-bounds} for precise definitions and a formal derivation of the bounds \eqref{eqn:transient-bounds} and \eqref{eqn:stationary-bounds}. 

We focus on bounding the stationary sensitivity in the context of stochastic reaction networks, i.e., $X(t, c)$ is a continuous time jump Markov process. 
To make use of the bounds \eqref{eqn:stationary-bounds}, we need reliable estimators for the IAF $\tau_{\pi^c}(f)$ and the path space FIM $\mathcal{I}_\mathcal{H}(P^{c})$.
For the IAF, we have shown in the Appendix \ref{app:UQ-bounds} that 
\[\tau_{\pi^c}(f) = \lim_{T\to \infty} \frac{1}{T}\text{Var}_{P_{[0,T]}^c}\left(\int_0^T f(X(s))\,ds\right).\]
Hence, when $T$ is large, an approximate estimator for the IAF is
\[\frac{1}{T(N-1)} \sum_{k=1}^N (Y^{(k)} - \bar{Y})^2,\]
where $N$ is the sample size, $Y^{(k)} = \int_0^T f(X^{(k)}(s)) ds$ is the $k$-th sample and $\bar{Y} = N^{-1}\sum_{k=1}^N Y^{(k)}$ is the sample average. 
Note that \eqref{eqn:stationary-bounds} assumes the dynamics starts at the stationary regime, hence a burn-in period is 
necessary for the dynamics to relax to the stationary state before we start sampling the IAF. 
Now for the path space FIM, it can be written as the stationary expectation of a special observable in terms of the propensity functions (see Appendix \ref{app:information-theory}),  i.e.,
\[\mathcal{I}_\mathcal{H}(P^{c}) = \mathbb{E}_{\pi^c}\left\{\sum_{j =1}^m \lambda_j(x, c)
\nabla \lambda_j(x, c) \nabla \lambda_j(x, c)^{\text{tr}} \right\}.\] 
Since the expectation is taken under the stationary distribution, the path space FIM can be 
approximated as the ergodic average of the observable. 
That is, 
\begin{equation*}
\begin{split}
&\mathcal{I}_\mathcal{H}(P^{c}) = \\
&\lim_{T\to\infty}\frac{1}{T}\int_0^T \sum_{j =1}^m \lambda_j(X(s), c)\nabla \lambda_j(X(s), c) \nabla \lambda_j(X(s), c)^{\text{tr}} \,ds.
\end{split}
\end{equation*}
Hence, an estimator for the path space FIM is simply
\[\frac{1}{N}\sum_{k=1}^N Z^{(k)},\]
where $Z^{(k)}$ is the $k$-th realization of the ergodic average.  
Note that the FIM is of great interests by itself since it reflects the identifiability of parameters by Cram\'er-Rao's inequality.  
We will use the path space information bounds $\eqref{eqn:stationary-bounds}$ to estimate
the stationary sensitivity bounds for numerical experiments in the next section.
\section{Numerical examples}\label{sec:numerical_eg}
In this section, we consider two bistable examples arising in chemistry and systems biology. 
We demonstrate 
that the ParRep algorithm can efficiently sample rare transitions between two stable equilibrium points and outperforms the standard SSA by a significant speedup factor.

\subsection{Bistable Schl\"{o}gl model}
\subsubsection{Model}
The Schl\"{o}gl model is one of the simplest example of stochastic reaction networks that exhibit bistability.
It is an auto-catalytic network involving three species whose population can change according to 
the reaction network in Table~\ref{tab:Schlogl}. 
Following our notational convention, we denote by $X_V(t)$ the concentration 
of the species $S$ and $X^V(t)$ the population of $S$.
The concentration of $A$ and $B$ (denoted by $a$ and $b$, respectively) are fixed due to an exchange of chemicals between two 
material baths \cite{vellela2009stochastic} and hence $a$ and $b$ are considered as parameters 
of the network. 
Therefore, it is equivalent to the Schl\"{o}glmodel as a one species network
\[2S \xrightleftharpoons[c_2]{c_1} 3S, \quad \varnothing \xrightleftharpoons[c_4]{c_3} S\] 
with $a$ and $b$ absorbed in the rate constants $c_1$ and $c_3$.
In this paper, we follow the chemical convention to write the reactions of Schl\"{o}gl network as in Table~\ref{tab:Schlogl}. 

\begin{table}[ht]
\label{tab:Schlogl}
\centering
\begin{tabular}{lll}
\hline
Reaction& Propensity Function   & Stoich. Vec. \\
\hline
\hline
$A + 2S \to 3S$ & $\lambda_1^V(x, c) = c_1 a x(x-1)/V$  & $\eta_1 = 1$\\
\hline
$3S \to A + 2S$  & $\lambda_2^V(x, c) = c_2 x(x-1)(x-2)/V^2$  & $\eta_2 = -1$\\
\hline
$B \to S$  & $\lambda_3^V(x, c) = c_3 b V$  & $\eta_3 = 1$\\
\hline
$S \to B$  & $\lambda_4^V(x, c) = c_4 x$  & $\eta_4 = -1$\\
\hline
\end{tabular}
\caption{Bistable Schl\"{o}gl model.}
\end{table}

\begin{figure}
\centering
\includegraphics[scale=0.5]{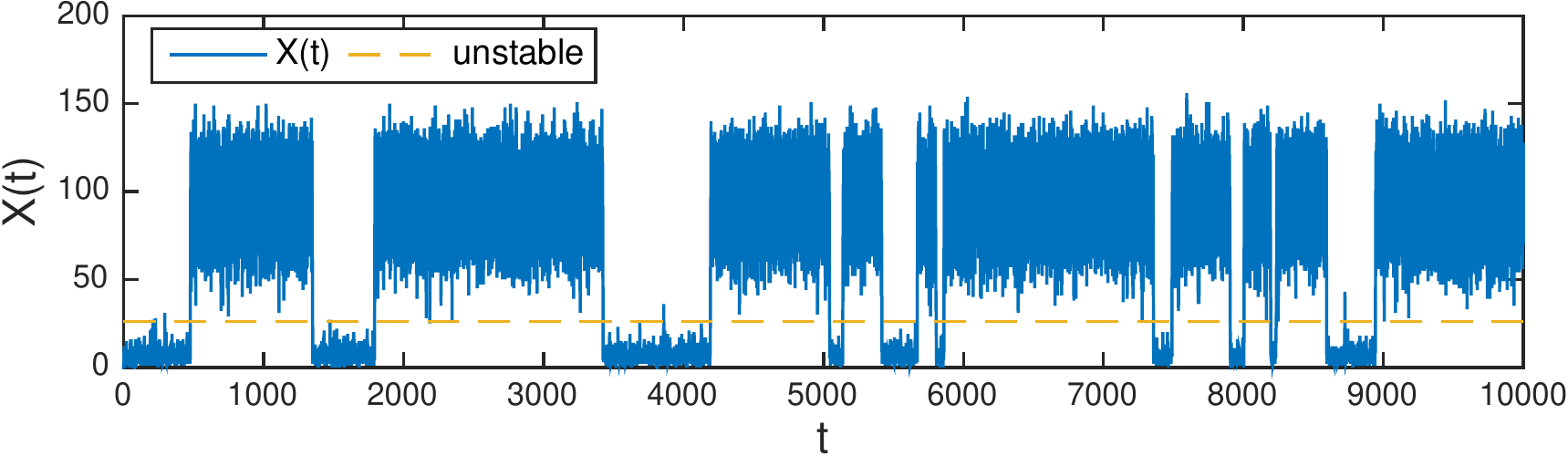}
\caption{A typical SSA trajectory of the Schl\"{o}gl model up to $t = 10^4$ with $V = 25$. The dashed line corresponds to the unstable equilibrium point multiplied by the system size $V = 25$, i.e., $\bar{X}^V = 25.9649$. Crossing of the unstable state $\bar{X}^V$ is a rare event.}
\label{fig:Schlogl-traject}
\end{figure}

In the large volume limit, the concentration process has a deterministic limit $\bar{x}$ satisfying the RRE \eqref{eqn:RRE}
\[\frac{d\bar{x}}{dt} = c_1 a \bar{x}^2 - c_2 \bar{x}^3 + c_3 b - c_4 \bar{x}.\]
We choose $a = 1$, $b = 2$, $c_1 = 3$, $c_2 = 0.6$, $c_3 = 0.25$ and $c_4 = 2.95$\cite{cao2013adaptively},
in which case the RRE has two stable equilibrium 
points $\bar{x}_+$ and $\bar{x}_-$ separated by an unstable equilibrium point $\bar{x}_0$. 
Therefore, Schl\"{o}gl model exhibits two time scales: the fast time scale corresponds to
the relaxation to one of the stable equilibrium points and the slow time scale corresponds to 
the rare transitions between the two stable equilibrium points.   
The two-time scale feature is illustrated in Figure~\ref{fig:Schlogl-traject}, where the standard SSA is performed with $V = 25$.

Due to the bistable nature, long time simulation is needed to sample enough transition events 
so that the system relaxes to stationary distribution. 
We apply the ParRep algorithm to accelerate the sampling of very long trajectories in order to estimate the stationary distribution $\pi^c$.
We decompose the state space into two metastable sets separated by the unstable equilibrium state $\bar{X}_0^V = V \bar{x}_0$
(we multiply the concentration $\bar{x}_0$ by $V$ so that all the comparisons are in terms of the population instead of the concentration). 
That is,
\[W_+ = \{x \in E: x \leq \bar{X}_0^V\}, \quad W_- =   \{x \in E: x > \bar{X}_0^V\},\]
where $E$ is the state space of $X^V(t)$. 
Note that this decomposition will be optimal for ParRep in terms of efficiency since both $W_+$ and $W_-$ will be strongly metastable. 
This can be seen by contradiction. In fact, if the decomposition is defined in terms of a point
$X'$ which is left to $V \bar{x}_0$ ($X' < V \bar{x}_0$), then every time $X^V$ exits from $W_+$ (with first exit state in the interval $(X', V \bar{x}_0)$) will be quickly pulled back to the left stable point $V\bar{x}_+$ with a dominating probability, by the large deviation principle. 
Hence, the subinterval $(X', V \bar{x}_0)$ in $W_-$ is not metastable and the ParRep will be inefficient since the parallel step is not activated when the process is in this interval.
Therefore, the optimal choices for separatrix is the point $V \bar{x}_0$ which guarantees that both of the decomposed sets are truly metastable.
\subsubsection{Results and discussion}

\begin{figure}
\centering
\includegraphics[scale=0.56]{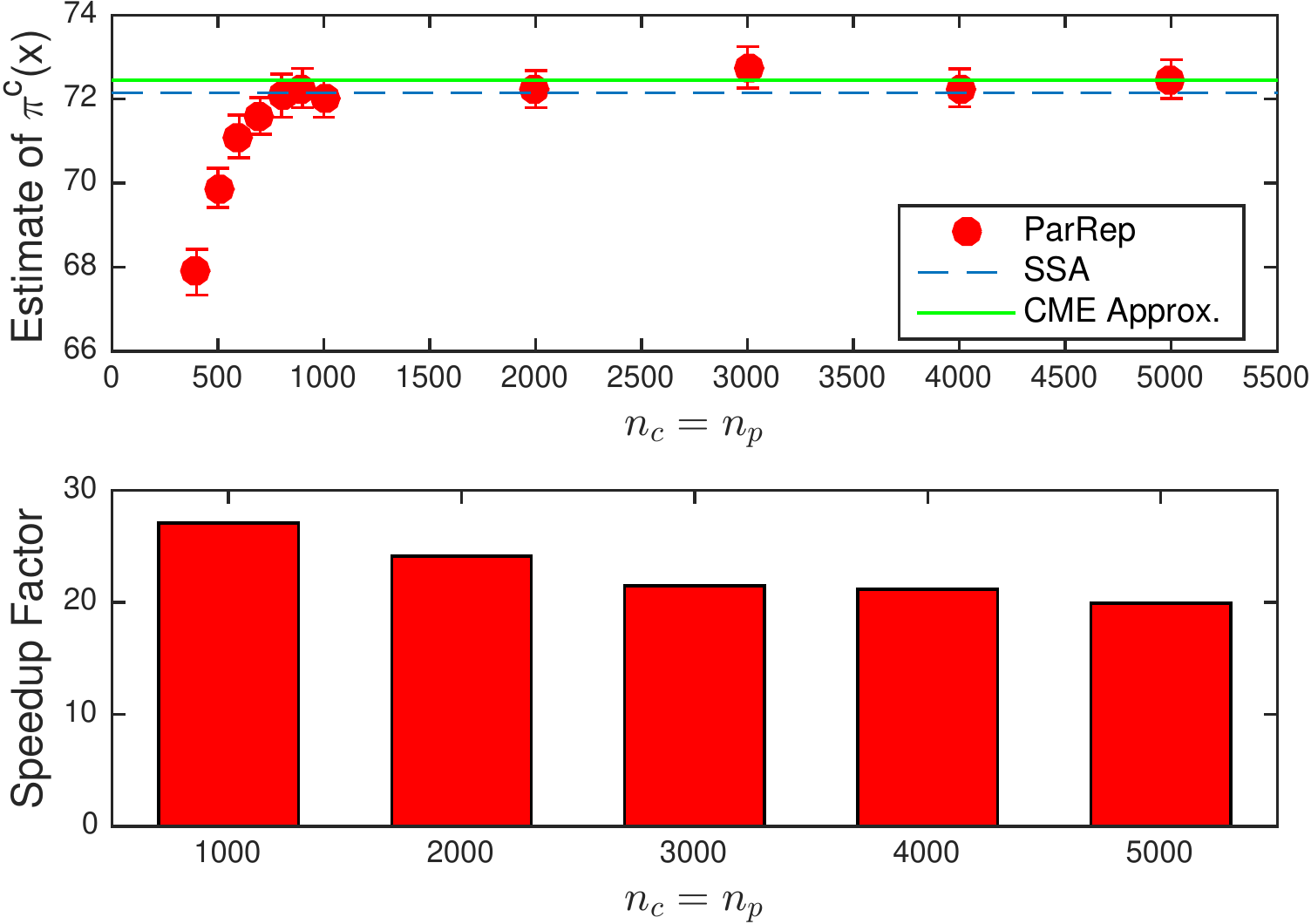}
\caption{Stationary average estimation of $\pi^c(x)$ (upper panel) and speedup factor (lower panel) with various decorrelation and dephasing steps. 
The speedup factor is only shown for $n_c = n_p = 1000, 2000, \ldots, 5000$.}
\label{fig:Schlogl-t-version}
\end{figure}

Figure~\ref{fig:Schlogl-t-version} shows the estimates of the stationary average of $X$  (ergodic average at $t = 10^5$) with SSA (blue dashed line) and with the ParRep algorithm (red dot with $95\%$ confidence interval) for different choice of decorrelation and dephasing steps.
The number of replicas for ParRep is $R =100$.
We also plot the numerical approximation of CME as a benchmark (green solid line) for accuracy in Figure~\ref{fig:Schlogl-t-version}.
It can be seen that the ParRep simulation approximates the stationary average very well (relative error with respect to the CME solution is $0.04\%$ for $n_c=n_p=5000$) when 
the decorrelation and dephasing steps are large, this is consistent with our expectation 
that the QSD of each metastable set is well approximated for large $n_c$ and $n_p$.
All simulation results are obtained based on $100$ sample trajectories.   
The CPU time of standard SSA simulation is about $192$ hours for $100$ samples. 
We demonstrate the corresponding speedup factor with  $n_c = n_p =1000$ to $5000$ (smaller values are ignored as the corresponding estimates are not accurate enough).
We can see that with $100$ replicas, our ParRep outperforms the standard SSA by a significant speedup factor.

\begin{figure}
\centering
\includegraphics[scale=0.56]{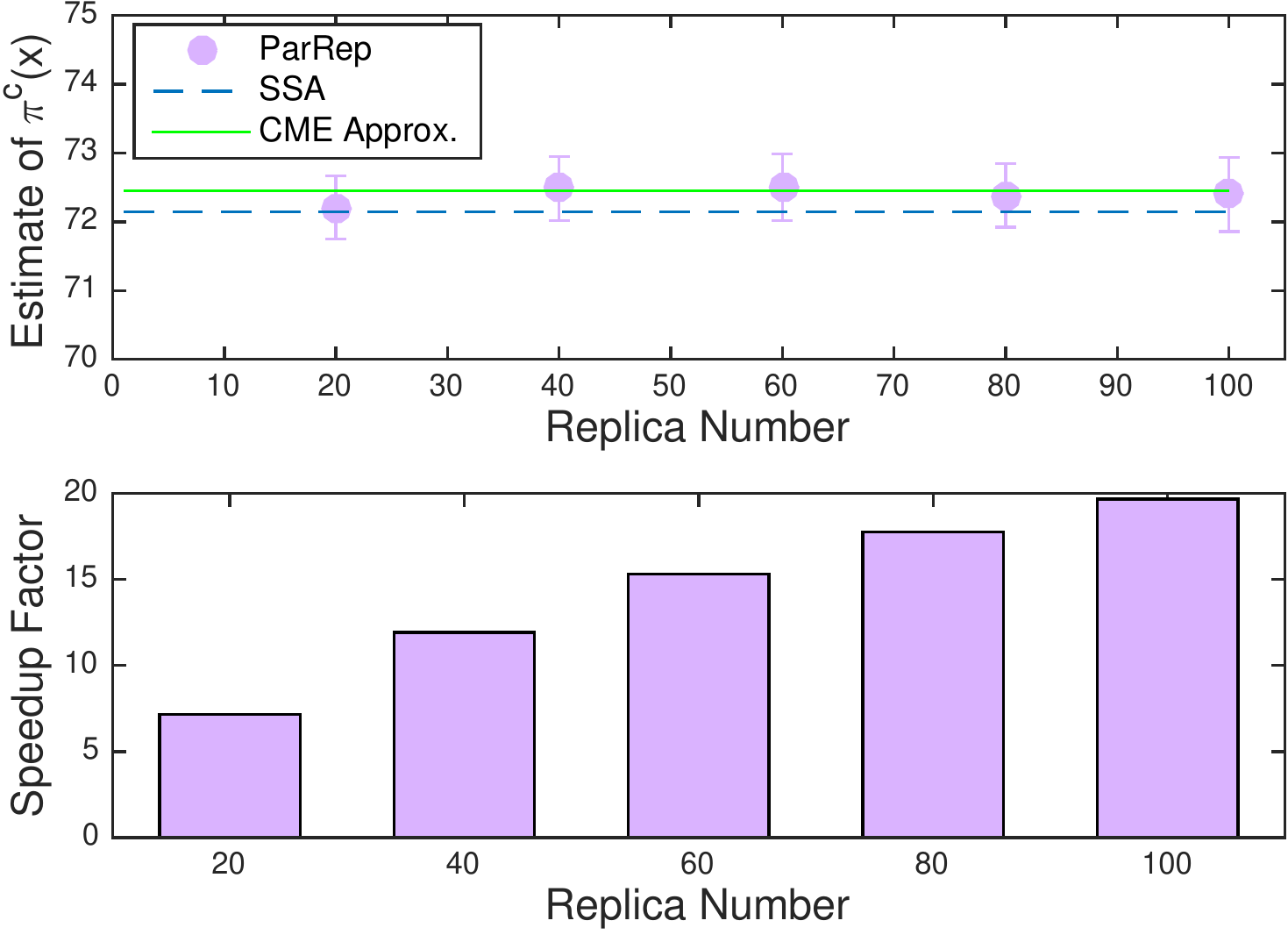}
\caption{Stationary average estimation of $\pi^c(x)$ (upper panel) and speedup factor (lower panel) with various number of replicas.}
\label{fig:Schlogl-replica-version}
\end{figure}

\begin{figure}
\centering
\includegraphics[scale=0.5]{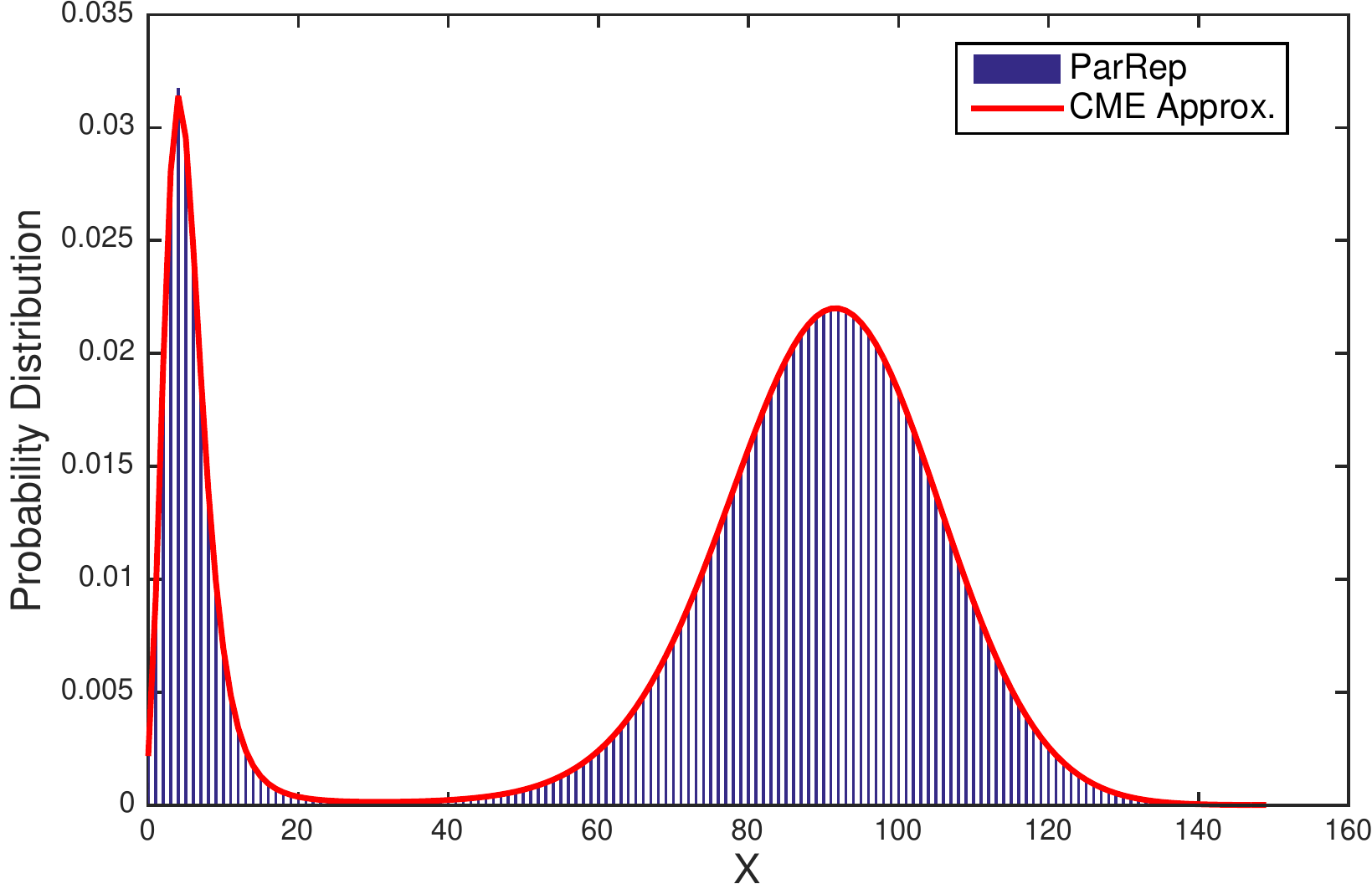}
\caption{Stationary probability distribution of $X^V$ with $V = 25$.
A total of $150$ bins (blue bar) are used to sample the empirical distribution on the interval $[0, 149]$.
The result is compared with the numerical solution (red curve) of the CME.  
}
\label{fig:Schlogl-distribution}
\end{figure}

We also study the efficiency of ParRep in terms of the number of replicas. 
In Figure~ \ref{fig:Schlogl-replica-version} we show the estimation of the stationary average of $X^V$ and the corresponding 
speedup factor. 
The decorreation and dephasing steps are fixed at $n_c = n_p = 5000$. 
We observe that the speedup factor increases from $7$ to $20$ when the number of 
replicas changes from $20$ to $100$.
However, the accuracy of ParRep is independent of the number of replicas.
In Figure~\ref{fig:Schlogl-distribution} we demonstrate the application of ParRep to estimate the probability distribution of $X$ with $n_c = n_p = 5000$. 
The estimated probability distribution (blue bar) is compared with the probability distribution
obtained from CME approximation.
The plot suggests that ParRep is a rather accurate method when suitable $n_c$ and $n_p$ are chosen.
  
Finally, we apply the path space information bound \eqref{eqn:stationary-bounds} to obtain a bound for the sensitivity index $S_{f, v}({\pi}^c)$. 
Here we only consider the stationary sensitivity of the observable $f(X^V(t)) = X^V(t)$ with respect to each parameter $c_i$, $i = 1, \ldots, 4$. 
Note that the stochastic reaction networks we simulate start at the transient regime, i.e., 
the initial distribution of $X^V(t)$ is not necessarily $\pi^c$. 
However, the path space information bounds \eqref{eqn:stationary-bounds} assume that $X^V(t)$ starts in 
the stationary regime. 
Therefore, a burn-in period is needed for the process to be well into the stationary regime
before we can start sampling the IAF $\tau_{\pi^c}(f)$ and the path space FIM $\mathcal{I}_{\mathcal{H}}(P^c)$.
We choose the terminal time $T = 2\times 10^5$ and use the first half $[0, 10^5]$ as
the burn-in period to prepare the stationary distribution and the second half $[10^5, 2\times 10^5]$ to sample the IAF and path space FIM.
The computed path space FIM and the confidence intervals are shown in Table \ref{tab:Schlogl_FIM}.
Note that the pFIM is not only useful for obtaining the final sensitivity bounds, but also implies the identifiability of parameters by the Cram\'er-Rao bound.
The computed IAF $\tau_{\pi^c}(f)$ is $5.87E$+$05$.
The resulting sensitivity bounds are shown in Table \ref{tab:Schlogl_SB}.
To see whether the obtained sensitivity bounds are tight enough, we compare them with the approximated sensitivities.
The approximation is obtained by differentiating the CME \eqref{eqn:CME} at steady state and truncating the state space to $[0, 149]$. 
The resulting equation is a linear system that can be solved numerically.
Comparing the sensitivity bounds with the approximated sensitivities, we observe that the bounds are not tight enough in this example.    
In fact, it has been observed in several examples that the path space information bounds are not always tight when applied to multi-scale problems.
Nevertheless, the bounds are quite useful for screening insensitive parameters in large scale stochastic dynamical systems.
We will demonstrate this application of the bounds in the next example. 

\begin{table}[ht]
\centering
\setlength{\tabcolsep}{8pt}
\begin{tabular}{c|c|c}
\hline
Matrix Element & Estimated pFIM & Half width C.I.  \\
\hline
\hline
$(1, 1)$ & 8.75E+01 & 3.02E-01 \\
\hline
$(2, 2)$ & 1.67E+03 & 5.66E+00 \\
\hline
$(3, 3)$ & 2.00E+02 & 2.59E-06 \\
\hline
$(4, 4)$ & 2.46E+01 & 7.88E-02\\
\hline
\end{tabular}
\caption{Estimated path space FIM for Schl\"{o}gl model}
\label{tab:Schlogl_FIM}
\end{table}
\begin{table}[ht]
\centering
\begin{tabular}{c|c|c|c|c}
\hline
Parameter & $c_1$ & $c_2$ & $c_3$ & $c_4$\\
\hline
\hline
Bounds & 7.16E+03  &  3.09E+04  & 1.08E+04  & 3.80+E03\\
\hline
CME Approx. & 4.07E+02  &  9.10E+02  & 6.30E+02  & -2.65+E02\\
\hline
\end{tabular}
\caption{Estimated sensitivity bounds and approximated sensitivities for Schl\"ogl model}
\label{tab:Schlogl_SB}
\end{table}


\subsection{Genetic switch with positive feedback}
\subsubsection{Model}
\begin{figure}[ht]
\begin{tikzpicture}
\node at (1,2.5) {$\text{DNA}_{\text{in}}$};
\node at (1,1) {$\text{DNA}_{\text{act}}$};
\node at (3.5,1) {$\text{mRNA}$};
\node at (6,1) {$\text{Protein}$};
\node at (3.5,2.5) {$\varnothing$};
\node at (6,2.5) {$\varnothing$};

\draw [->] (0.9,2.3) -- (0.9,1.3);  
\draw [->] (1.1,1.3) -- (1.1,2.3);  
\draw [->] (1.6,1) -- (2.9,1); 
\draw [->] (4.1,1) -- (5.4, 1); 
\draw [->] (3.5,1.3) -- (3.5,2.3); 
\draw [->] (6,1.3) -- (6, 2.3); 

\node at (0.3, 1.7) {$F(X_2)$};
\node at (1.8, 1.7) {$G(X_2)$};
\node at (2.2,1.2) {$a$};
\node at (4.7,1.2) {$\gamma b$};
\node at (3.2,1.7) {$\gamma$};
\node at (5.8,1.7) {$1$};
\end{tikzpicture}
\caption{Diagram of the genetic switch network. The switching rates $F$ and $G$ are
given in \eqref{eqn:on-off-rates}.}
\label{fig:genetic-switch-diagram}
\end{figure}
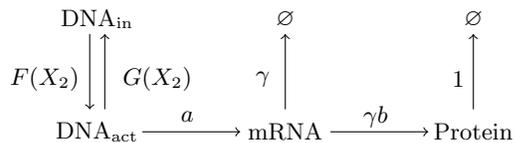

Another example we study in this paper is the genetic switch network 
which is the fundamental mechanism for cells to shift between alternate 
gene-expression states.
See Figure~\ref{fig:genetic-switch-diagram} for the diagram of the network.
In the genetic switch network, there is an on-off switch
for DNA to be in the active or inactive state.
Hence the total population of active DNA and inactive DNA is $1$.
The transition rates $F$ (inactive to active) and $G$ (active to inactive) between these two states depend on the 
number of proteins through  a positive feedback.
Following Assaf, Roberts and Luthey-Schulten (2012)\cite{assaf2011determining},
we explicitly take the mRNA noise into account
since it has been shown that the presence of mRNA has a significant impact on 
the dynamics of the network \cite{mehta2008exponential}.
We list the propensity function and stochiometric vector of each reaction channel in Table \ref{tab:genetic-switch}.


\begin{table*}[ht]
\centering
\setlength{\tabcolsep}{12pt}
\begin{tabular}{lll}
\hline
Reaction& Propensity Function & Stoich. Vec. \\
\hline
\hline
$\text{DNA}_{\text{act}} \to \text{DNA}_{\text{act}} + \text{mRNA}$ & $\lambda_1^V(\xi, x_1, x_2) = a\xi$  & $\eta_1 = (1, 0) $\\
$\text{mRNA} \to \varnothing$ & $\lambda_2^V(\xi, x_1, x_2) = \gamma x_1$ & $\eta_2 = (-1, 0) $  \\
$\text{mRNA} \to \text{mRNA}+\text{Protein}$ & $\lambda_3^V(\xi, x_1, x_2) = \gamma b x_1$ & $\eta_3 = (0, 1) $  \\
$\text{Protein} \to \varnothing$ & $\lambda_4^V(\xi, x_1, x_2) = x_2$  & $\eta_4 = (0, -1) $\\
\hline
\end{tabular}
\caption{Genetic switching system}
\label{tab:genetic-switch}
\end{table*}

We fix large volume $V = ab = 2400$ throughout this example.
The two-dimensional process $X^V(t) = (X_1^V(t), X_2^V(t))$ denotes the number of mRNA and protein at time $t$.
We denote the number of active DNA by the process $\xi(t)$ and hence the number of inactive DNA is $1-\xi(t)$.
The transition rates are taken to be of the Hill-type functions $F(x_2)$ for the 
inactive to active transition and $G(x_2)$ for the reverse transition, where
\begin{equation}\label{eqn:on-off-rates}
\begin{split}
F(x_2) &= k_0^{\text{min}} + (k_0^{\text{max}} - k_0^{\text{min}}) x_2^2/(x_2^2 + D^2)\\
G(x_2) &= k_1^{\text{max}} - (k_1^{\text{max}} - k_1^{\text{min}}) x_2^2/(x_2^2 + D^2).
\end{split}
\end{equation}
Throughout this example, we follow Assaf et al. \cite{assaf2011determining} to set the parameters as follows: $a = 2400/b$, $b = 22.5$, $\gamma = 50$,  $k_0^{\text{min}} = k_1^{\text{min}} = 24/b$, $k_0^{\text{max}} = k_1^{\text{max}} = 2400/b$ and $D = 1000$.

We point out that the genetic switch model does not fall into the standard framework of stochastic reaction networks we describe in Sec \ref{sec:SRN}.
In fact, the random time change representation of $X^V(t)$ can be written as
\begin{equation}
X^V(t) = X^V(0) + \sum_{j=1}^{4}\mathcal{P}_j\left(\int_0^t \lambda_j^V(\xi(s), X^V(s))\,ds\right) \eta_j,
\end{equation}
See Table \ref{tab:genetic-switch} for the four reactions involved in this representation. 
Note that the propensity functions $\lambda_j$ are functions of both the switching variable 
$\xi(t)$ and the population process of mRNA and protein $X^V(t)$.
Since $\xi(t) \in \{0, 1\}$ does not depend on volume $V$,
the process $(\xi, X^V(t))$ does not satisfy the large volume limit \eqref{eqn:fluid-limit}.
However, the mean numbers of mRNA and protein still satisfy the following rescaled RRE (i.e., the ODE governing $\bar{X}^V(t) = V\bar{x}(t)$) \cite{lv2014constructing, li2015large}
\begin{equation}\label{eqn:genetic-switch-ODE}
\begin{split}
\frac{d\bar{X}_1}{dt}& = \frac{a F(\bar{X}_2)}{F(\bar{X}_2) + G(\bar{X}_2)} - \gamma \bar{X}_1\\
\frac{d\bar{X}_2}{dt} &= \gamma b \bar{X}_1 - \bar{X}_2,
\end{split}
\end{equation}
where the factor $F/(F+G)$ gives the probability that the DNA is in an active state. 
With our choice of parameters, \eqref{eqn:genetic-switch-ODE} has two stable equilibrium 
points $\bar{X}_+^V = (0.0225, 25. 2628)$ and $\bar{X}_-^V = (1.6353, 1839.6883)$ separated by a saddle point $\bar{X}_0^V = (0.4545, 511.2865)$.
Therefore, the genetic switching network is bistable. When $V$ is finite, there are noise induced rare transitions between $\bar{X}_+^V$ and $\bar{X}_-^V$.

To find the optimal decomposition of the state space $E$ into two metastable sets, we need to 
analyze the phase space of \eqref{eqn:genetic-switch-ODE} to determine the separatrix of the two metastable regions induced by $\bar{X}_+^V$ and $\bar{X}_-^V$.
Unlike the Schl\"{o}gl model , it is nontrivial to find the separatrix in this example since it is in 
$\mathbb{R}^2$.
Instead, the way we detect rare transitions is ad-hoc.  
We simply choose the horizontal line that passes through the saddle point $\bar{X}_0^V$
as the boundary to define the two metastable sets. 
We provide a heuristic explanation for the choice. 
From the large deviation perspective, there exists a most probable transition path from 
$\bar{X}_+^V$ to $\bar{X}_-^V$ \cite{li2015large,lv2014constructing,dykman1994large} such that if a transition occurs, then with a dominating probability, the transition would move along this path.
We know that the true separatrix passes through the saddle point $\bar{X}_0^V$ and that the most probable transition path crosses the true separatrix along a path that is ``sufficiently close" to the saddle point (\cite{lv2014constructing}).
Since the points $(0, 511)$, $(0, 512)$, $(1, 511)$ and $(1, 512)$ are the only possible states that are sufficiently close to the saddle point $(0.4545, 511.2865)$, the most probable transition path can only cross the separatrix (from $\bar{X}_+^V$ to $\bar{X}_-^V$) by moving from $(0, 511)$ to
$(0, 512)$ or to $(1, 511)$ depending on where the true separatrix lies. 
Accordingly, this suggests that we can use either $y = 511. 2865$ (if the most probable path moves from $(0, 511)$ to $(0, 512)$) or $x = 0.4545$ (if the most probable path moves from $(0, 511)$ to $(1, 511)$) as the boundary to decompose the state space into two metastable sets. 
It is readily seen that we should choose the horizontal one since
the process is much more metastable in the $y$ direction than that in the $x$ direction.
Therefore, we decompose the state space into two sets
\[W_+ = \{(x, y) \in E| y < 511.2865\}\]
and 
\[W_- = \{(x, y) \in E| y > 511.2865\}.\]
Our simulation results confirm that this is a good choice of decomposition.
Note that though the choice of decomposition may not be optimal since we do not know the true separatrix a priori, it only affects the efficiency but not the accuracy of ParRep as we discuss in Section \ref{sec:ParRep-accuracy-efficiency}.
A rigorous approach to defining the optimal  decomposition into metastable sets is the subject of ongoing work.

\begin{figure}[ht]
\centering
\includegraphics[scale=0.6]{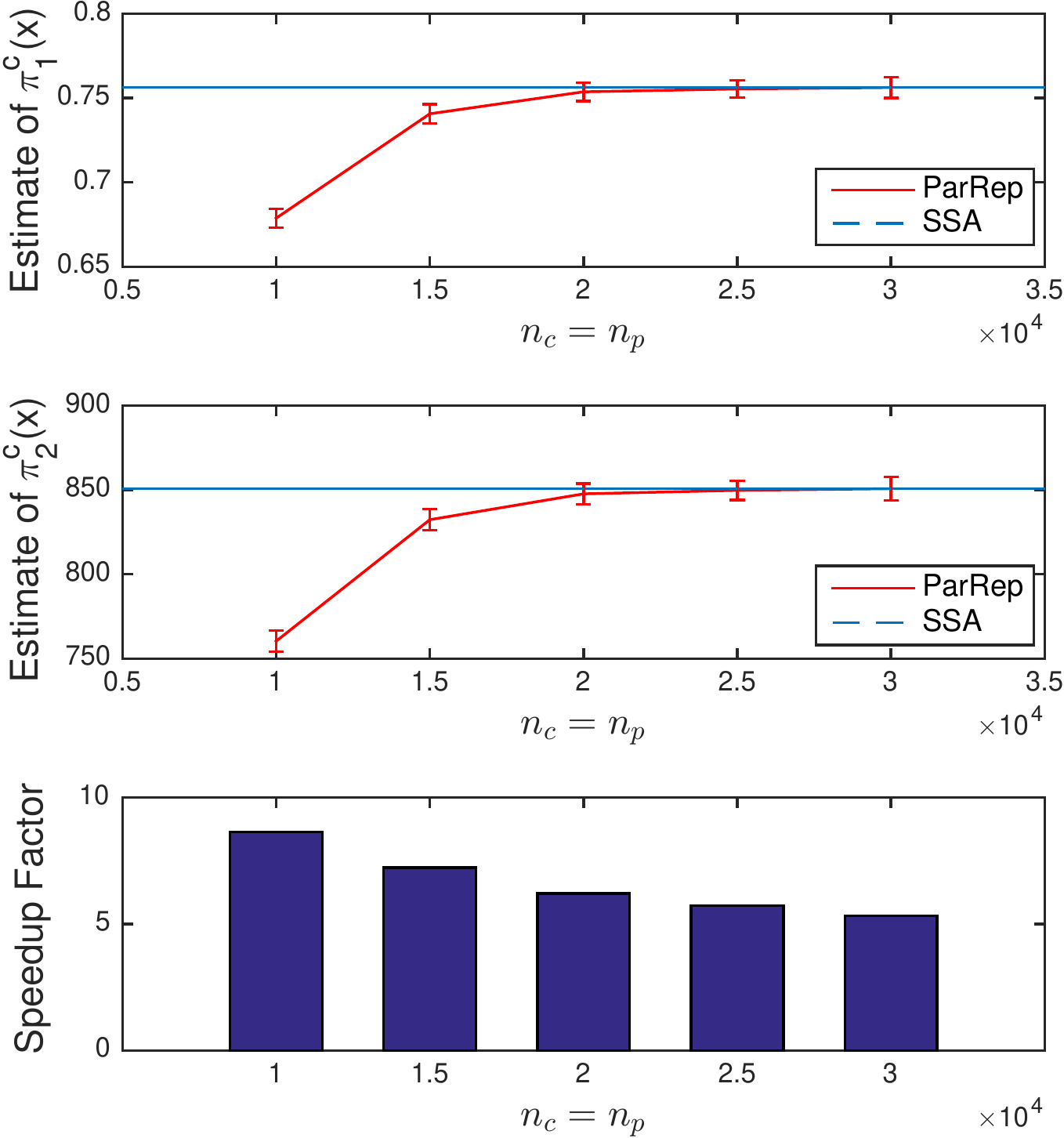}
\caption{Stationary average estimations of $\pi_1^c(x)$ and $\pi_2^c(x)$  (upper panel) and speedup factor (lower panel) with various number of replicas, where $\pi_1^c$ and $\pi_2^c$ are the marginal distribution with respect to $x_1$ and $x_2$, respectively.}
\label{fig:genetic-switch-t-version}
\end{figure}

\begin{figure}[ht]
\centering
\includegraphics[scale=0.6]{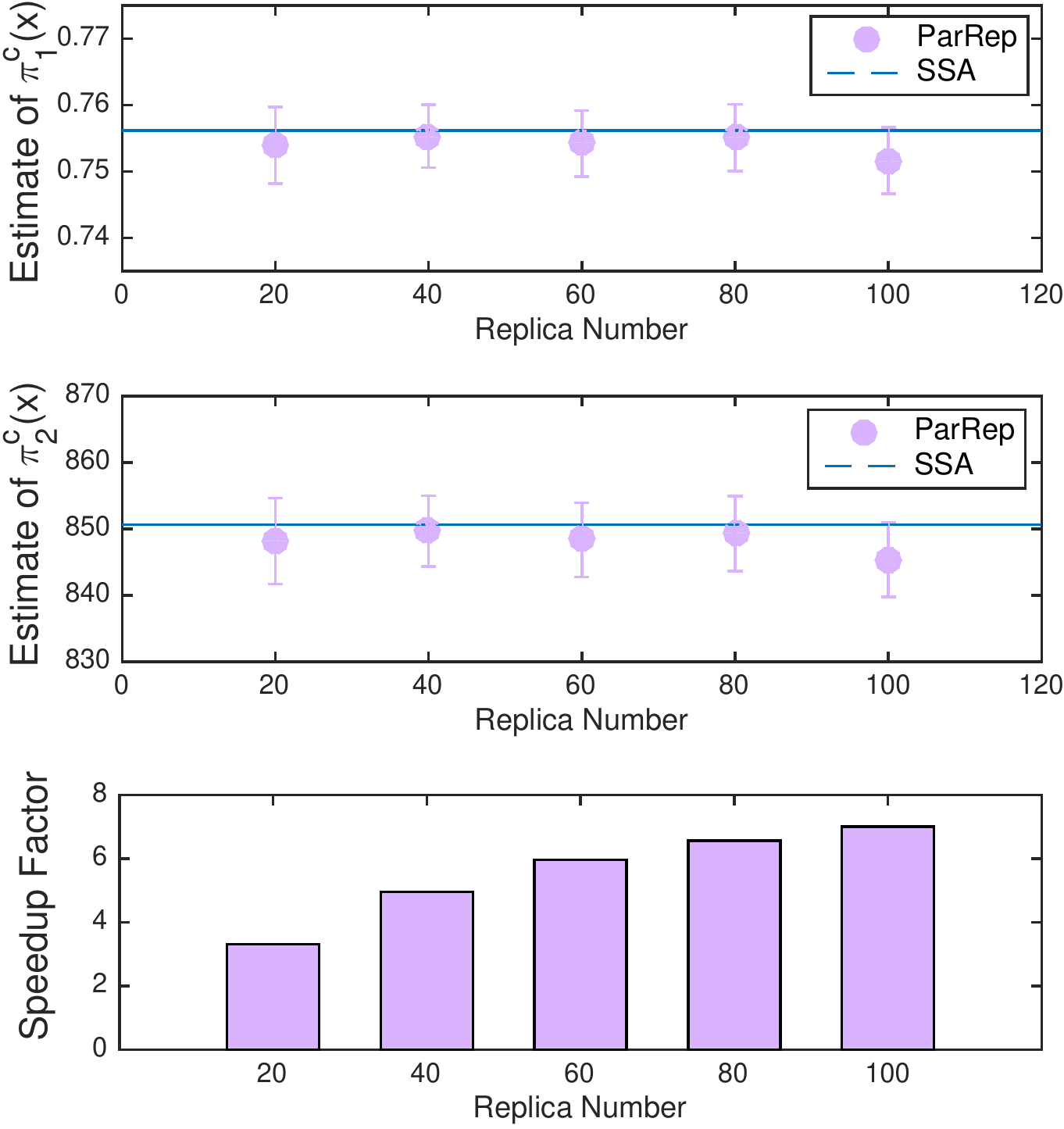}
\caption{Stationary average estimation of $\pi_1^c(x)$ and $\pi_2^c(x)$ (upper panel) and speedup factor (lower panel) with various number of replicas, where $\pi_1^c$ and $\pi_2^c$ are the marginal distribution with respect to $x_1$ and $x_2$, respectively.}
\label{fig:genetic-switch-replica-version}
\end{figure}

\subsubsection{Results and discussion}
Throughout the simulation of the genetic switch network, all simulation results are obtained based on $100$ sample trajectories.
The initial population is $1$ molecule for inactive DNA and $0$ molecule for all the remaining  species.
The terminal time is taken to be $T = 10^6$ which is sufficiently large  for sampling the ergodic average. 
We first study the accuracy of ParRep in terms of the decorrelation step $n_c$ and dephasing step $n_p$ with $100$ replicas. Figure \ref{fig:genetic-switch-t-version} demonstrates the simulation results regarding the stationary means of mRNA and protein when $n_c = n_p$ increase. 
Simulation results with SSA (blue dashed line) are used for accuracy comparison.  
The corresponding speedup factor is shown in the same plot (lower panel). 
The plot suggests that with $n_c=n_p = 2\times10^4$ or above, the ParRep is as accurate
as the standard SSA. 
Figure \ref{fig:genetic-switch-replica-version} shows the speedup of ParRep when we vary the number of replicas with $n_c=n_p = 2\times10^4$. 
We do not gain as much speedup as in the Schl\"{o}gl model since the genetic switch network requires
 much larger $n_c$ and $n_p$ to converge to the QSD at each metastable set (as we observed in Figure~ \ref{fig:genetic-switch-t-version}).
Nevertheless, we can see that with $100$ replicas, the speedup factor of ParRep is about $7\times$ when compared to SSA.

\begin{figure}[ht]
\centering
\includegraphics[scale=0.45]{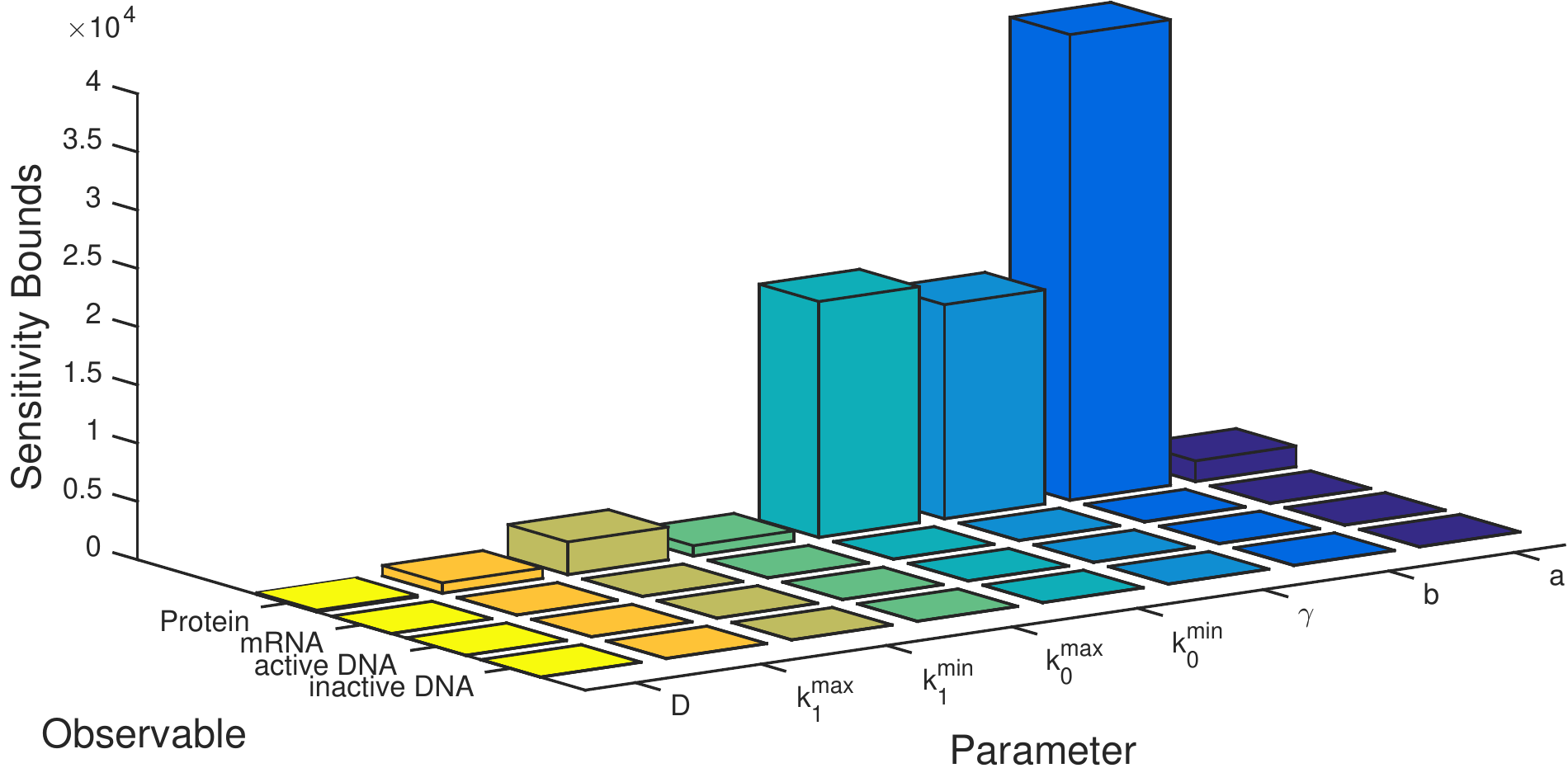}
\caption{Path space sensitivity bounds of genetic switches network.}
\label{fig:genetic-switch-SB}
\end{figure}

We also study the parametric sensitivity of the genetic switch model by quantifying the 
stationary path space sensitivity bounds \eqref{eqn:stationary-bounds}.
The observables in considerations are the number of active DNA ($\xi$), inactive DNA ($1-\xi$), mRNA ($X_1^V$) and protein ($X_2^V$). 
The parameters are arranged in the order 
$a, b, \gamma, k_0^{\text{min}}, k_0^{\text{max}}, k_1^{\text{min}}, k_1^{\text{max}}, D$.
We aim to apply ParRep to estimate the stationary sensitivity bounds of each observable with respect to each of the parameters. 
In order to obtain the bounds, we simulate the process up to final time $2\times10^6$. 
The time interval $[0, 10^6]$ corresponds to the transient regime and $[10^6, 2\times10^6]$ 
corresponds to the stationary regime.
The estimated (diagonal) path space FIM along with the confidence interval are shown in Table \ref{tab:genetic-switch-pFIM}.
The estimated IAF for each observable is shown in Table \ref{tab:genetic-switch-IAT}. 
Finally, we combine the path space FIM and IAF to obtain the stationary sensitivity bounds. 
To illustrate the observation that most of the sensitivity indices are small, we visualize the sensitivity bounds in Figure \ref{fig:genetic-switch-SB}. 
We see that the the active DNA, inactive DNA and mRNA are insensitive to parametric perturbation, whereas the protein tends to be sensitive. 
If we are interested in quantifying the parametric uncertainty for the genetic switch model, these
sensitivity bounds suggest that we can screen out those insensitive combinations and apply direct methods to estimate the remaining sensitivity such as the number of protein with respect to $b$, 
$\gamma$ and $k_0^{\text{min}}$.
Note that without this bounds, we have to estimate $4 \times 8 = 32$ sensitivities even when we do not take other observables into consideration. 
However, with the sensitivity bounds for screening, we only need to estimate much fewer sensitivities depending on the controlled confidence level we use.
Therefore, this two-step strategy significantly reduces the computational cost especially when applied to large scale networks.


\begin{table}[ht]
\centering
\setlength{\tabcolsep}{8pt}
\begin{tabular}{c|c|c}
\hline
Matrix Element & Estimated pFIM & Half width C.I. \\
\hline
\hline
$(1, 1)$  & 3.34E-03    &    2.35E-05  \\
\hline
$(2, 2)$  & 1.69E+00   & 1.19E-02\\
\hline
$(3, 3)$  & 3.57E-01    &   2.52E-03\\
\hline
$(4, 4)$  & 4.34E-01   &    2.83E-03\\
\hline
$(5, 5)$  & 8.32E-04    &  5.88E-06\\
\hline
$(6, 6)$  & 8.33E-03   &  6.12E-05\\
\hline
$(7, 7)$  & 8.44E-04   &  5.23E-06\\
\hline
$(8, 8)$  & 2.22E-05  &  1.56E-07\\
\hline
\end{tabular}
\caption{Path space FIM of genetic switch network.}
\label{tab:genetic-switch-pFIM}
\end{table}
\begin{table}[ht]
\centering
\setlength{\tabcolsep}{6pt}
\begin{tabular}{cccc}
\hline
active DNA & inactive DNA & mRNA  & Protein\\
\hline
\hline
1.64E+02 &1.64E+02 & 7.47E+02 & 9.45E+08\\
\hline
\end{tabular}
\caption{IAFs of genetic switch network.}
\label{tab:genetic-switch-IAT}
\end{table}


\appendix 
\section{Basics for information theory}\label{app:information-theory}
We review some basic information theory concepts for the completeness of the paper.
In particular, we briefly reproduce the formula for the relative entropy and the path space FIM in the context of stochastic reaction networks (i.e., continuous time jump Markov processes) \cite{pantazis2013relative}.

Given the path space probability distribution $P_{[0, T]}^c$ and its perturbation 
$P_{[0, T]}^{c'}$ on the path space $E$, 
their pseudo-distance can be measured by the 
relative entropy
\begin{equation*}
\mathcal{R}(P_{[0, T]}^c  | P_{[0, T]}^{c'}) = \int_E \log \left(\frac{dP_{[0, T]}^c(x)}{dP_{[0, T]}^{c'}(x)}\right)  P_{[0, T]}^c(dx).
\end{equation*}
In particular, the Radon-Nikodym derivative follows from the following Girsanov formula \cite{bremaud1981point}
\begin{equation*}
\begin{split}
\frac{dP_{[0, T]}^c(x)}{dP_{[0, T]}^{c'}(x)} = & \frac{\mu^c(x(0))}{\mu^{c'}(x(0))} \prod_{j=1}^m \exp\left(\int_0^T \log\frac{\lambda_j(x(s),c)}{\lambda_j(x(s),c')}dR_j(s)\right. \\
&\left.- \int_0^T \lambda_j(x(s), c)  - \lambda_j(x(s),c') ds\right),
\end{split}
\end{equation*}
where $R_j(t)$ is the count of the $j$-th reaction up to time $t$.
Assuming the dynamics starts from the stationary regime (i.e., $\mu^c = \pi^c$) and using the facts that $R_j(t) - \int_0^t \lambda_j(X(s)) ds$ is a martingale under $P_{[0, T]}^c$, 
the relative entropy can be simplify as
\[\mathcal{R}(P_{[0, T]}^c  | P_{[0, T]}^{c'}) = T\mathcal{H}(P^c  | P^{c'}) + \mathcal{R}(\mu^c | \mu^{c'}), \]
where 
\begin{equation*}
\begin{split}
&\mathcal{H}(P^c  | P^{c'}) = \\
&\mathbb{E}_{\mu^c}\left\{\sum_{j=1}^m \lambda_j(x,c) \log \frac{\lambda_j(x,c)}{\lambda_j(x,c')} - \lambda_j(x,c) + \lambda_j(x, c') \right\}
\end{split}
\end{equation*}
is the path space relative entropy rate (RER). Note that when $c' = c+ \delta$, the Taylor expansion of
$\mathcal{H}(P^c  | P^{c'})$ gives 
\begin{equation*}
\begin{split}
\mathcal{H}(P^c  | P^{c'}) 
=\frac{1}{2} \delta^{\text{tr}} \mathcal{I}_\mathcal{H}(P^c)\delta + \mathcal{O}(|\delta|^3)
\end{split}
\end{equation*}
where 
\[\mathcal{I}_\mathcal{H}(P^c)=\mathbb{E}_{\mu^c}\left\{\lambda_j(x, c) \nabla_c \log \lambda_j(x, c) \nabla_c \log \lambda_j(x, c)^{\text{tr}}\right\}\]
is the path space Fisher information matrix (FIM) of RER $\mathcal{H}(P^c  | P^{c'}) $.

\section{Path space information bounds: a formal derivation}\label{app:UQ-bounds}
For completeness of the paper, we give a formal derivation of the path space information bounds
on both the transient regime and the stationary regime, see the reference \cite{dupuis2016path} for rigorous derivation.  
We consider a continuous time Markov process $X(t, c)$ with stationary distribution $\pi^c$. 
For the path space measure $P_{[0,T]}^c$, we assume that it is absolutely continuous with respect to a reference measure $R_{[0, T]}$ such that
$P_{[0,T]}^c(dx) = p_{[0, T]}^c(x){R_{[0, T]}(dx)}$ for any $c$.
Then by the definition of sensitivity indices and the Cauchy-Schwarz inequality,  
\begin{equation*}
\begin{split}
&S_{f, v}\left(P_{[0, T]}^c\right) 
= \left|\lim_{\epsilon \to 0} \frac{1}{\epsilon}\left(\mathbb{E}_{P_{[0,T]}^{c+\epsilon v}}[f] - \mathbb{E}_{P_{[0,T]}^c}[f]\right)\right| \\
=& \left|\lim_{\epsilon \to 0}\frac{1}{\epsilon}\int_E f(x) \left(\frac{p_{[0,T]}^{c+\epsilon v}(x)}{p_{[0,T]}^c(x)} - 1\right) P_{[0, T]}^c(dx)\right|\\
=& \left|\int_E \left(f(x) - \mathbb{E}_{P_{[0,T]}^c}[f] \right)\frac{\left.\frac{d}{d\epsilon}\right|_{\epsilon = 0}p_{[0,T]}^{c+\epsilon v}(x)}{p_{[0,T]}^c(x)} P_{[0, T]}^c(dx)\right|\\
\leq & \sqrt{\text{Var}_{P_{[0,T]}^{c}}(f)} \sqrt{v^{\text{tr}} \mathcal{I}(P_{[0, T]}^c) v},
\end{split}
\end{equation*}
where $\mathcal{I}(P_{[0, T]}^c)$ is the FIM of the relative entropy $\mathcal{R}(P_{[0,T]}^{c}|P_{[0,T]}^{c+\epsilon v})$.
This gives the path space information bounds on the transient regime.

In the stationary regime, we focus on the ergodic average type observable $F(x) = T^{-1}\int_0^T f(x(s)) \,ds$. 
Since the stationary distribution $\pi^c$ is also the initial distribution of the stochastic process $X(t, c)$, it holds that $S_{F, v}(P_{[0, T]}^c) = S_{f, v}(\pi^c)$.
Hence by the path space information bounds for $F$,  
\begin{equation*}
\begin{split}
S_{f, v}(\pi^c) &\leq \sqrt{\text{Var}_{P_{[0,T]}^{c}}(F)} \sqrt{v^{\text{tr}} \mathcal{I}(P_{[0, T]}^c) v}\\
&=  \sqrt{\frac{1}{T}\text{Var}_{P_{[0,T]}^{c}}(TF)} \sqrt{\frac{1}{T}v^{\text{tr}} \mathcal{I}(P_{[0, T]}^c) v}.
\end{split}
\end{equation*}
Taking $T \to \infty$, 
\begin{equation*}
\begin{split}
S_{f, v}(\pi^c) &\leq  \sqrt{\lim_{T\to \infty}\frac{1}{T}\text{Var}_{P_{[0,T]}^{c}}(TF)} \sqrt{v^{\text{tr}} \mathcal{I}_{\mathcal{H}}(P^c) v},
\end{split}
\end{equation*}
where we assumed that $T^{-1}\mathcal{I}(P_{[0, T]}^c)$ converges to the path space FIM $\mathcal{I}_{\mathcal{H}}(P^c)$. 
Note that 
\begin{equation*}
\begin{split}
&\lim_{T\to \infty}\frac{1}{T}\text{Var}_{P_{[0,T]}^{c}}(TF) \\
=&\lim_{T\to \infty} \frac{1}{T}\mathbb{E}_{P_{[0, T]}^c}\left\{\left(\int_0^T f(x(s)) - \pi^c(f)\, ds\right)^2\right\}\\
= &\lim_{T\to \infty}\frac{1}{T}\int_0^T \int_0^T \text{Cov}_f(u-v, 0) \,du\,dv\\
=&\int_{-\infty}^{\infty} \text{Cov}_f(s, 0) \,ds
\end{split}
\end{equation*}
gives the integrated auto-correlation function (IAF) for observable $f$, where
\[\text{Cov}_f(u, v) = \mathbb{E}_{P_{[0, T]}^c}\left\{(f(x(u)) - \pi^c(f))(f(x(v)) - \pi^c(f))\right\}\]
is the covariance between $f(x(u))$ and $f(x(v))$.
We denote that IAF by
\[\tau_{\pi^c}(f)=\lim_{T \to \infty}\frac{1}{T}\text{Var}_{P_{[0,T]}^{c}}(TF).\]
We remark that the IAF only differs with the integrated auto-correlation time (IAT) by a multiplying factor $\text{Cov}_f(0, 0)$, i.e., $\text{IAF} = \text{Cov}_f(0, 0)\times\text{IAT}$. 
Finally, we have the stationary path space information bounds
\[S_{f, v}(\pi^c) \leq  \sqrt{\tau_{\pi^c}(f)} \sqrt{v^{\text{tr}} \mathcal{I}_{\mathcal{H}}(P^c) v}.\]

\begin{acknowledgements}
The work of P.P. has been partially supported by the U.S. Department of Energy, Office of Science, 
Office of Advanced Scientific Computing Research, Applied Mathematics program under the contract 
number DE-SC0010549. 
The work of T.W.  was partially supported by the DARPA project W911NF-15-2-0122.
We thank Professor Tiejun Li for discussions about the genetic switches example. 
\end{acknowledgements}

\bibliographystyle{abbrv}
\bibliography{ParRep_bistable}

\end{document}